\documentclass[twoside,11pt,reqno]{amsart}

\usepackage{amsmath,amsthm,amssymb,amstext,amsfonts,amscd}
\usepackage{graphicx,color}

\usepackage{multirow}

\numberwithin{equation}{section}
\begin{document}

\title[On a generalization of l'Hopital's rule]{On a generalization of l'Hopital's rule for multivariable functions}

\author[{\bf  V. V. Ivlev and I. A. Shilin}]{\bf V. V. Ivlev and I. A. Shilin$^*$}

\address{Department of Mathematics, Sholokhov Moscow State University for the Humanities,
Verhnya Radishevskaya 16-18, Moscow 109240, Russia}
\email{ilyashilin@li.ru}

\bigskip



\begin{abstract}
For students and their lecturers and instructors interested in the
natural problem of a possible generalization of l'Hopital's rule
for functions depending on two or more variables, we offer our
approach.  For instructors, we discuss the technique of
constructing indeterminate forms at a given point and having a
given double limit.
 \end{abstract}

 \thanks{$^*$ Corresponding author}

\maketitle

 \section{Introduction}\medskip

 In order to compute a multivariable function limit taking an indeterminate form, we use some well known ways.
 First, we can change Cartesian coordinate system in the plane to polar
 system $(x,y)\longmapsto(r,\varphi)\equiv\left(\sqrt{x^2+y^2},\arctan\frac{y}x\right)
 $: for example, $$\lim\limits_{(x,u)\to(0,0)}
\frac{x^2+y^2}{\sqrt{x^2+y^2+1}-1}=\left\{\frac00\right\}=
\lim\limits_{r\to0}\frac{r^2}{\sqrt{r^2+1}-1}\\=
\lim\limits\left(\sqrt{r^2+1}+1\right)=2. $$ In case of three or
more variables we use the spherical coordinate system. Second, we
use simply the definition of a limit: for example, the hypothesis
$\lim\limits_{(x,y)\to(0,0)}\frac{x^2y}{x^2+y^2}=0 $ is true
because, for arbitrary $\epsilon>0$, in the neighborhood
$\{(x,y)\,\mid\,0<\sqrt{x^2+y^2}<\epsilon\}$ we have
$$\left|\frac{x^2y}{x^2+y^2}-0\right|=\frac{x^2|y|}{x^2+y^2}\le
|y|\le\sqrt{x^2+y^2}<\epsilon.$$ Third, we can use some algebraic
transformations: for instance, $$\lim\limits_{(x,u)\to(2,2)}
\frac{x^2-y^2}{x^2-2xy+y^2+x-y}=\left\{\frac00\right\}=
\lim\limits_{(x,y)\to(2,2)}\frac{x+y}{x-y+1}=4 . $$ Fourth, we use
the estimation of function values: for instance, for each
$\epsilon$-neighborhood
$U(\infty,\epsilon):=\{(x,y)\,\mid\,\sqrt{x^2+y^2}>\epsilon^{-1}\}
$ of infinite point (the motivation of such definition is related
to the stereographic projection of the Riemann sphere to the
plane) we have $$\left| \frac{2x^2+2y^2+x+y}{3x^2+3y^2} \right|<
\frac23+\frac1{3|x|}+\frac1{3|y|}< \frac{2(1+\epsilon)}3 ,$$
therefore, $
\lim\limits_{(x,y)\to\infty}\frac{2x^2+2y^2+x+y}{3x^2+3y^2}
=\frac23 $.

In order to show that the limit does not exist, we compute the
associated directional limits: along the curve $\gamma_1$,
consisting of points  $(t,t^4) $, we have $$
\lim\limits_{(x,y)\to\infty}\frac{x^4+y^2}{x^2+y^4}=
\lim\limits_{t\to\infty}\frac{t^2+t^6}{t^{14}+1}=0,$$ along the
curve  $\gamma_2$, consisting of points  $(t^2,t) $, we have $$
\lim\limits_{(x,y)\to\infty}\frac{x^4+y^2}{x^2+y^4}=
\lim\limits_{t\to\infty}\frac{t^6+1}{2t^2}=+\infty.$$

But there is the question: how to show the existence of a limit?
In case of a single variable function, we often use the well known
l'Hopital's rule. Typically our students are interested in an
analogue of this rule with respect to multivariable functions.
Unfortunately, the existing textbooks do not contain any
corresponding theorems. In this paper, we suggest our
generalization of the group of theorem known as l'Hopital's rule.
For readers which are interested more deeply in this topic, we
recommend the paper \cite{l} which contains another approach to
the problem discussed here.

\section{A generalization of l'Hopital's rule}\medskip

The function $f(x,y) $  is called infinitely small of the first
order at point $(x_0,y_0) $,  if
$\lim\limits_{(x,y)\to(x_0,y_0)}f(x,y)=0 $  and
$[f'_x(x_0,y_0)]^2+[f'_y(x_0,y_0)]^2\ne0 $  (the last condition
simply means that ${\rm d}f(x_0,y_0)\ne0 $).

Let $n$  be a positive integer and  $n>1$. The function $f(x,y) $
is called infinitely small of order $n$  at point $(x_0,y_0) $ if
$\lim\limits_{(x,y)\to(x_0,y_0)}f(x,y)=0 $,   ${\rm
d}^if(x_0,y_0)=0 $ for  $i\in\{1,\ldots,n-1\} $, and ${\rm
d}^nf(x_0,y_0)\ne0 $.

Let us call the fraction $\frac{f(x,y)}{g(x,y)} $  as
indeterminate form of order $n$ at point $(x_0,y_0) $  if  $f$ and
$g$ are infinitely small of order $n$  at this point.

\textbf{Theorem 1.} \quad \emph{Let $f(x,y)$  and  $g(x,y)$ be two
functions defined, differentiable, and infinitely small of the
first order in some neighborhood of the point  $(x_0,y_0)$. Let,
in addition, $f'_x(x_0,y_0)\ne0 $ and  $f'_y(x_0,y_0)\ne0 $. The
non-zero limit $f'_x(x_0,y_0)\ne0 $ exists if and only if
\begin{equation}\label{f1}\frac{f'_x(x_0,y_0)}{g'_x(x_0,y_0)}=
\frac{f'_y(x_0,y_0)}{g'_y(x_0,y_0)}=k.\end{equation}}

\begin{proof} {\sl Necessity}. Let us introduce the following
 denotations: $k_1:=
\frac{f'_x(x_0,y_0)}{g'_x(x_0,y_0)}
  $  and  $k_2:=
\frac{f'_y(x_0,y_0)}{g'_y(x_0,y_0)}
  $. Since  $f$ and  $g$ are infinitely small at $(x_0,y_0)$, we can write the corresponding
  Taylor series in the above neighborhood of the point $(x_0,y_0)$
  as \begin{gather*}
f(x_0+\Delta x,y_0+\Delta y)= f'_x(x_0,y_0)\Delta x+
f'_y(x_0,y_0)\Delta y+o_1\left(\sqrt{(\Delta x)^2+(\Delta
y)^2}\right),
\\
g(x_0+\Delta x,y_0+\Delta y)= g'_x(x_0,y_0)\Delta x+
g'_y(x_0,y_0)\Delta y+o_2\left(\sqrt{(\Delta x)^2+(\Delta
y)^2}\right).
\end{gather*} Then \begin{multline*} k= \lim\limits_{(x,y)\to(x_0,y_0)} \frac{f(x,y)}
{g(x,y)}\\=
\lim\limits_{(x,y)\to(x_0,y_0)}\frac{f'_x(x_0,y_0)\Delta x+
f'_y(x_0,y_0)\Delta y}{g'_x(x_0,y_0)\Delta x+ g'_y(x_0,y_0)\Delta
y}\\=\lim\limits_{(x,y)\to(x_0,y_0)}\frac{k_1g'_x(x_0,y_0)\Delta
x+ k_2g'_y(x_0,y_0)\Delta y}{g'_x(x_0,y_0)\Delta x+
g'_y(x_0,y_0)\Delta y}.
 \end{multline*} In case of $\Delta y=o(\Delta x) $  we have $\lim\limits_{\Delta x\to0}\frac{\Delta y}{\Delta x}=0 $  and, therefore,  $k=k_1 $.
  The opposite case $\Delta x=o(\Delta y) $  gives us  $k=k_2 $.
 It means that the limit  $\lim\limits_{(x,y)\to(x_0,y_0)}\frac{f(x,y)}
{g(x,y)} $ does not exist.

{\sl Sufficiency}. If the condition  $
\frac{f'_x(x_0,y_0)}{g'_x(x_0,y_0)}=
\frac{f'_y(x_0,y_0)}{g'_y(x_0,y_0)}=k
  $ holds true, then $$ \lim\limits_{(x,y)\to(x_0,y_0)}\frac{f(x,y)}
{g(x,y)}=k
\lim\limits_{(x,y)\to(x_0,y_0)}\frac{g'_x(x_0,y_0)\Delta x+
g'_y(x_0,y_0)\Delta y}{g'_x(x_0,y_0)\Delta x+ g'_y(x_0,y_0)\Delta
y}=k.$$ The proof is complete. \end{proof}

For example, let us consider the limit
$\lim\limits_{(x,y)\to(1,1)}\frac{x^2+2xy-3y^2}{x^3-y^3} $. Here
we have the indeterminate form  $\left\{\frac00\right\} $. All the
first order partial derivatives of $f:=x^2+2xy-3y^2 $ and
$g:=x^3-y^3 $  at point $(1,1)$  are finite and do not equal to
zero:
\begin{gather*}f'_x(1,1)=2x+2y\mid_{(1,1)}=4,\,\,\,\,\,\,\,\,
f'_y(1,1)=2x-6y\mid_{(1,1)}=-4,\\
g'_x(1,1)=3x^2\mid_{(1,1)}=3,\,\,\,\,\,\,\,\,\,\,\,
g'_y(1,1)=-3y^2\mid_{(1,1)}=-3.\end{gather*} According to Theorem
1,
$$\lim\limits_{(x,y)\to(1,1)}\frac{x^2+2xy-3y^2}{x^3-y^3}=
\frac{f'_x(1,1)}{g'_x(1,1)}= \frac{f'_y(1,1)}{g'_y(1,1)}=\frac43.
 $$

 Let us note that (\ref{f1}) provides the existence the associated repeated limits. Moreover, in this
 case $$ \lim\limits_{x\to x_0} \lim\limits_{y\to y_0}
\frac{f(x,y)}{g(x,y)}= \lim\limits_{y\to y_0} \lim\limits_{x\to
x_0} \frac{f(x,y)}{g(x,y)}= \lim\limits_{(x,y)\to (x_0,y_0)}
\frac{f(x,y)}{g(x,y)}.
  $$ However, Theorem 1 does not cover all possible cases related to indeterminate forms. The following questions remain open:

  \begin{enumerate}

  \item Only one of the derivatives  $f'_x(x_0,y_0) $,  $f'_y(x_0,y_0) $,  $g'_x(x_0,y_0) $,  $g'_y(x_0,y_0) $ is equal to zero.

  \item $f'_x(x_0,y_0)=g'_x(x_0,y_0)=0 $ or $f'_y(x_0,y_0)=g'_y(x_0,y_0)=0
  $.

  \item $f'_x(x_0,y_0)=g'_y(x_0,y_0)=0 $ or $f'_y(x_0,y_0)=g'_x(x_0,y_0)=0
  $.

  \item $f'_x(x_0,y_0)=f'_y(x_0,y_0)=0 $ or $g'_x(x_0,y_0)=g'_y(x_0,y_0)=0
  $.

 \item  All partial derivatives of the first order are equal to zero.

 \end{enumerate}

 Let us consider successively these possibilities.

 (1) Without loss of generality, assume that  $f'_x(x_0,y_0)=0 $. We
 have \begin{multline*}
\lim\limits_{(x,y)\to(x_0,y_0)}\frac{f(x,y)} {g(x,y)}=
\lim\limits_{(x,y)\to(x_0,y_0)}\frac{ f'_y(x_0,y_0)\Delta y }
{g'_x(x_0,y_0)\Delta x+g'_y(x_0,y_0)\Delta y}\\=
\lim\limits_{(x,y)\to(x_0,y_0)}\frac{ f'_y(x_0,y_0) }
{g'_x(x_0,y_0)\frac{\Delta x}{\Delta y}+g'_y(x_0,y_0)}.
 \end{multline*} Let $\Delta x\to0 $  and $\Delta y\to0 $  such that  $\frac{\Delta x}{\Delta y}=r $.
  If $r=0$  then  $\lim\limits_{(x,y)\to(x_0,y_0)}\frac{f(x,y)}
{g(x,y)}=\frac{ f'_y(x_0,y_0)}{g'_y(x_0,y_o)} $. In case of
$r=\infty$, we have  $\lim\limits_{(x,y)\to(x_0,y_0)}\frac{f(x,y)}
{g(x,y)}=0 $. It means that the double limit does not exist, but
the repeated limits are exist and differ.

(2) In case of  $f'_x(x_0,y_0)=g'_x(x_0,y_0)=0 $, we have $$
\lim\limits_{(x,y)\to(x_0,y_0)}\frac{f(x,y)} {g(x,y)}=
\lim\limits_{(x,y)\to(x_0,y_0)}\frac{ f'_y(x_0,y_0)\Delta y }
{g'_y(x_0,y_0)\Delta y}= \frac{ f'_y(x_0,y_0) } {g'_y(x_0,y_0)}.
$$ Thus the limit exists if $
\lim\limits_{(x,y)\to(x_0,y_0)}\frac{ f'_x(x_0,y_0)\Delta x }
{g'_x(x_0,y_0)\Delta x}= \frac{ f'_y(x_0,y_0) } {g'_y(x_0,y_0)}. $

(3) In case of  $f'_x(x_0,y_0)=g'_y(x_0,y_0)=0 $, the double limit
does not exist because of $$
\lim\limits_{(x,y)\to(x_0,y_0)}\frac{f(x,y)} {g(x,y)}=
\lim\limits_{(x,y)\to(x_0,y_0)}\frac{ f'_y(x_0,y_0)\Delta y }
{g'_x(x_0,y_0)\Delta x},
$$ the repeated limits differ: $$
\lim\limits_{\Delta x\to0}\lim\limits_{\Delta y\to0}f(x,y)=0
,\,\,\,\,\,\, \lim\limits_{\Delta y\to0}\lim\limits_{\Delta
x\to0}f(x,y)=\infty.
 $$

 (4) The case $f'_x(x_0,y_0)=f'_y(x_0,y_0)=0 $  means that at point $(x_0,y_0) $
   the function $f$  is infinitely small of order $\ge2$  and the function $g$  is infinitely small of the first order.
    It follows to $
\lim\limits_{(x,y)\to(x_0,y_0)}\frac{f(x,y)}
{g(x,y)}=\lim\limits_{\Delta x\to0}\lim\limits_{\Delta
y\to0}f(x,y) =\lim\limits_{\Delta y\to0}\lim\limits_{\Delta
x\to0}f(x,y)=0.
 $

 (5) Let  $
f'_x(x_0,y_0)=f'_y(x_0,y_0)= g'_x(x_0,y_0)=g'_y(x_0,y_0)= 0 $.
Under the circumstances we deal with an indeterminate form of the
second order.

\textbf{Theorem 2.} \quad \emph{Let  $f(x,y)$ and $g(x,y)$  be two
functions defined, twice differentiable, and infinitely small of
the second order in some neighborhood of the point  $(x_0,y_0)$.
Let, in addition, that all their partial derivatives of the second
order are not equal to zero. The non-zero limit
$k:=\lim\limits_{(x,y)\to(x_0,y_0)}\frac{f(x,y)} {g(x,y)} $ exists
if and only if
\begin{equation}\label{f2}\frac
{ f"_{xx}(x_0,y_0) }{g"_{xx}(x_0,y_0)} =\frac { f"_{xy}(x_0,y_0)
}{g"_{xy}(x_0,y_0)} =\frac { f"_{yx}(x_0,y_0) }{g"_{yx}(x_0,y_0)}
=\frac { f"_{yy}(x_0,y_0) }{g"_{yy}(x_0,y_0)}=k.
\end{equation}}

\begin{proof} For convenience, we introduce the following
matrices: $$
(a_{ij}):=\begin{pmatrix}f"_{xx}(x_0,y_0)& f"_{xy}(x_0,y_0)\\
f"_{yx}(x_0,y_0)&f"_{yy}(x_0,y_0)\end{pmatrix},
(b_{ij}):=\begin{pmatrix}g"_{xx}(x_0,y_0)& g"_{xy}(x_0,y_0)\\
g"_{yx}(x_0,y_0)&g"_{yy}(x_0,y_0)\end{pmatrix}.
 $$ It is well known that a double limit does not depend on the
direction. Let us consider the direction  $\gamma_r$ such that
$\Delta y=r\Delta x $. If (\ref{f2}) holds true, the double limit
is equal to $k$  along all lines $\gamma_r$ containing point
$(x_0,y_0)$: \begin{multline*}
\lim\limits_{(x,y)\to(x_0,y_0)}\frac{f(x,y)} {g(x,y)}\\=
\lim\limits_{(x,y)\to(x_0,y_0)} \frac{a_{11}(\Delta x)^2+
[a_{12}+a_{21}]\Delta x\Delta y+ a_{22}(\Delta y)^2+o_1}
{b_{11}(\Delta x)^2+ [b_{12}+b_{21}]\Delta x\Delta y+
b_{22}(\Delta y)^2+o_2}\\
=\lim\limits_{(x,y)\to(x_0,y_0)}\frac {k(b_{11}+ [b_{12}+b_{21}]r+
b_{22}r^2)} {b_{11}+ [b_{12}+b_{21}]r+ b_{22}r^2}=k,
 \end{multline*} where $o_1$  and $o_2$  denote the corresponding remainder
 terms of Taylor's formula. In opposite, if at least two of the fractions included in (\ref{f2}) differ, we have a limit depending on
 $r$:\begin{multline*}
\lim\limits_{(x,y)\to(x_0,y_0)}\frac{f(x,y)} {g(x,y)}\\=
\lim\limits_{(x,y)\to(x_0,y_0)} \frac{a_{11}(\Delta x)^2+
[a_{12}+a_{21})]\Delta x\Delta y+ a_{22}(\Delta y)^2+o_1}
{b_{11}(\Delta x)^2+ [b_{12}+b_{21}]\Delta x\Delta y+
b_{22}(\Delta y)^2+o_2}\\
=\lim\limits_{(x,y)\to(x_0,y_0)}\frac {a_{11}+ [a_{12}+b_{21}]r+
a_{22}r^2} {b_{11}+ [b_{12}+b_{21}]r+ b_{22}r^2}.
 \end{multline*} The proof is complete. \end{proof}

 For example, the limit $\lim\limits_{(x,y)\to(0,0)}\frac{x^2+xy+y^2}
{x^2-xy+y^2}$  does not exist because the numerator and
denominator are infinitely small functions of the first order and
$ \frac{f"_{xx}(0,0)}{g"_{xx}(0,0)}=1=-
\frac{f"_{xy}(0,0)}{g"_{xy}(0,0)} $.

 It is clear that Theorems 1 and 2 could be generalized for the indeterminate forms $\left\{\frac00\right\} $  of the  $n$-th order:

\textbf{Theorem 3.} \quad \emph{Let  $f(x,y)$ and $g(x,y)$  be two
functions defined, $n$ times differentiable, and infinitely small
of the second order in some neighborhood of the point $(x_0,y_0)$.
Let, in addition, that all their partial derivatives of order $n$
are not equal to zero. The non-zero limit
$k:=\lim\limits_{(x,y)\to(x_0,y_0)}\frac{f(x,y)} {g(x,y)} $ exists
if and only if
\begin{equation}=\frac {\frac{\partial^nf}{\partial
x^l\,\partial^{n-l}}\mid_{(x_0,y_0)}} {\frac{\partial^ng}{\partial
x^l\,\partial^{n-l}}\mid_{(x_0,y_0)}}=k\,\,\,\,\,\,\text{for all
$l\in\{0,1,\ldots,n\}$}
\end{equation} (we suppose here that all partial derivatives of the  $n$-th order are
continuous).}

From the above theorems, it is possible to derive the
corresponding corollaries for the indeterminate forms
$\left\{\frac\infty\infty\right\} $,  $\{0\cdot\infty\} $, etc.
For instance, in case of  $\left\{\frac\infty\infty\right\} $, we
have $$ \lim\limits_{(x,y)\to(x_0,y_0)}\frac{f(x,y)}{g(x,y)}=
\left\{\frac\infty\infty\right\}=
\lim\limits_{(x,y)\to(x_0,y_0)}\frac{[g(x,y)]^{-1}}
{[f(x,y)]^{-1}}=\left\{\frac00\right\}.
$$ Assuming that  $x_0,y_0\in\bar{\mathbb{R}}=\mathbb{R}\bigcup\{\pm\infty\} $,
we can obtain the corresponding corollary for limits  $
\lim\limits_{(x,y)\to(\infty,\infty)}\frac{f(x,y)} {g(x,y)} $:
changing the variables  $x$ and $y$  to $u:=x^{-1} $  and
$v:=y^{-1} $, we have $$
\lim\limits_{(x,y)\to(\infty,\infty)}\frac{f(x,y)}{g(x,y)}=
\lim\limits_{(x,y)\to(0,0)} \frac{f\left(\frac1u,\frac1v\right)}
{f\left(\frac1u,\frac1v\right)}.
$$ Finally, all considered cases could be generalized for functions of three and more variables.

\section{Construction of indeterminate forms for examples and tasks}\medskip

For an instructor, it is very important to construct a good
collection of appropriate examples and tasks providing the
successful learning of the above generalization of l'Hopital's
rule. In this section, we discuss how to do it. Let us note first
that if $f(x,y)$  is a $n$  times differentiable function in a
neighborhood of the point  $(x_0,y_0)$, then $f(x,y)$  can be
expressed as the Taylor series:
$$f(x,y)=f(x_0,y_0)+\sum\limits_{k=1}^n\,\frac{{\rm
d}^kf\mid_{(x_0,y_0)}}{k!}+r_n.$$  From this follows that $ \Delta
f\mid_{(x_0,y_0)} \equiv f(x,y)-f(x_0,y_0) $ is a infinitely small
function of the first order in case of  ${\rm
d}f\mid_{(x_0,y_0)}\ne0 $. Further,  $\Delta f\mid_{(x_0,y_0)}-
{\rm d}f\mid_{(x_0,y_0)} $ is a infinitely small function of the
second order in case of  ${\rm d}^2f\mid_{(x_0,y_0)}\ne0 $. And so
on. We obtain that for any functions $f(x,y)$  and $g(x,y)$, which
are differentiable in a neighborhood of the point $(x_0,y_0)$, the
limit
\begin{equation}\label{f3}\lim\limits_{(x,y)\to(x_0,y_0)}\frac{\Delta f\mid_{(x_0,y_0)}}{
\Delta g\mid_{(x_0,y_0)}}\end{equation} is an indeterminate form
of the first order. Similarly, for any twice differentiable
functions $f(x,y)$  and  $g(x,y)$  , the limit
\begin{equation}\label{f4}\lim\limits_{(x,y)\to(x_0,y_0)}\frac{\Delta f\mid_{(x_0,y_0)}
-{\rm d}f\mid_{(x_0,y_0)} }{ \Delta g\mid_{(x_0,y_0)}-{\rm
d}g\mid_{(x_0,y_0)} }
\end{equation} is an indeterminate form of the second order. And so on. However, it is not clear whether the limits
 (\ref{f3}) and (\ref{f4}) exist. Assume that fractions $\frac{f'_x(x_0,y_0)}{g'_x(x_0,y_0)} $  and $\frac{f'_y(x_0,y_0)}{g'_y(x_0,y_0)} $
   differ. According to Theorem 1, the condition \begin{equation}\label{f5}\frac{f'_x(x_0,y_0)}{g'_x(x_0,y_0)}
=\frac{f'_y(x_0,y_0)}{g'_y(x_0,y_0)} \end{equation} is necessary
and sufficient for the existence of finite limit (\ref{f1}) which
we denote as  $k$. Let us add the linear terms $C_1(x-x_0) $  and
$C_2(y-y_0) $ to  $f(x,y)$. Than (\ref{f3}) can be rewritten as  $
\frac{f'_x(x_0,y_0)+C_1}{g'_x(x_0,y_0)}
=\frac{f'_y(x_0,y_0)+C_2}{g'_y(x_0,y_0)} $. It is easy to see that
$C_1=kg'_x(x_0,y_0)-f'_x(x_0,y_0) $ and
$C_2=kg'_y(x_0,y_0)-f'_y(x_0,y_0) $. According to Theorem 2, the
condition $$ \frac{f"_{xx}(x_0,y_0)}{g'_{xx}(x_0,y_0)} =
\frac{f"_{xy}(x_0,y_0)}{g"_{xy}(x_0,y_0)}=
\frac{f"_{yx}(x_0,y_0)}{g"_{yx}(x_0,y_0)}=
\frac{f"_{yy}(x_0,y_0)}{g"_{yy}(x_0,y_0)}
 $$ is criterion for the existence of finite limit (\ref{f2}) which we denote again  as  $k$.
  Repeating the same approach, we rewrite $f(x,y)$  as
  $$f(x,y)+C_1^*(x-x_0)^2+C_2^*(x-x_0)(y-y_0)+C_3^*(y-y_0)^2.$$ It follows to
  formulas \begin{gather*} C_1^*=\frac{ k g"_{xx}(x_0,y_0)- f"_{xx}(x_0,y_0)}2,\,\,\,
\,\,\,\,
C_3^*=\frac{ k g"_{yy}(x_0,y_0)- f"_{yy}(x_0,y_0)}2,\\
C_2^*=k g"_{xy}(x_0,y_0)- f"_{xy}(x_0,y_0)
\end{gather*} (under the supposition that the mixed derivatives are continuous and, therefore, coincide).

For example, let $f:=x^2y+x+y $  and $g:=x^2y^2+xy $  and we need
to construct the indeterminate form  $\left\{\frac00\right\} $ of
the first order at point $(1,1)$  under the condition that the
corresponding limit is equal to  $2$. Both our functions are
infinitely differentiable, \begin{gather*} \Delta
f\mid_{(1,1)}=x^2y+x+y-3,\,\,\,\,\,
\Delta g\mid_{(1,1)}=x^2y^2+xy-2,\\
C_1=2(2xy^2+y)-2xy-1\mid_{(1,1)}=3,\,\,
C_2=2(2x^2y+x)-x^2-1\mid_{(1,1)}=4.
\end{gather*} So we can offer for students the limit  $\lim\limits_{(x,y)\to(1,1)}
\frac{x^2y+x+y+3(x-1)+4(y-1)-3}{x^2y^2-2} $, that is
$\lim\limits_{(x,y)\to(1,1)} \frac{x^2y+4x+5y-10}{x^2y^2-2} $, as
an example illustrating the above generalization of l'Hopital's
rule or as a task for classroom work or homework. In order to
construct the indeterminate form of the second order, we obtain
that \begin{gather*} \Delta f\mid_{(1,1)}-
{\rm d}f\mid_{(1,1)}=x^2y-2x-y+2,\\
\Delta g \mid_{(1,1)} -
{\rm d}g\mid_{(1,1)}=x^2y^2+xy-3x-3y+4,\\
C_1^*=2y^2-y\mid_{(1,1)} =1,\,\,\,\,
C_3^*=2x^2\mid_{(1,1)} =2,\\
C_2^*=8xy-2x+2\mid_{(1,1)} =8.
\end{gather*} So we have
 $$
\lim\limits_{(x,y)\to(1,1)}\frac{x^2y+x^2+8xy-12x+2y^2
-13y+13}{x^2y^2+xy-3x-3y+4}=2,
$$
where first order derivatives of the numerator and denominator are
equal to zero at point $(1,1)$.

\section{Conclusion}\medskip

The approach discussed in this paper we successfully use in our
teaching practice for math first year students at Sholokhov Moscow
State University for the Humanities. We hope, this approach could
be useful for all lectures, instructors and student teaching and
learning, relatively, mathematical analysis. Some of ideas
discussed in this paper, we used in the book \cite{i}

\bigskip

\renewcommand{\refname}{

\end{document}